\documentclass[12 pt]{article}
\newtheorem{Theorem}{Theorem}
\newtheorem{Lemma}{Lemma}
\newtheorem{Corollary}{Corollary}
\begin{document}
\title{Asymptotically normal estimators for Zipf's law}
\author{Mikhail Chebunin\thanks{E-mail: chebuninmikhail@gmail.com, Novosibirsk State University, Novosibirsk, Russia},
Artyom Kovalevskii \thanks{E-mail: kovalevskiii@gmail.com, Novosibirsk State Technical University, 
Novosibirsk State University, Novosibirsk State University of Economics and Management, Novosibirsk, Russia.  
The research was supported by RFBR grant 17-01-00683}}
\date{}
\maketitle
\begin{abstract}
Zipf's law states that sequential frequences of words in a text correspond to a power function. Its probabilistic model is an infinite urn scheme with
asymptotically power distribution. The exponent of this distribution must be estimated. We use the number of different words in a text and similar statistics to construct asymptotically normal estimators of the exponent.
\end{abstract}


\section{Introduction}

Zipf's law (Zipf, 1949) states that sequential frequences $f_i$  
of words in a text are equal to $ci^{-1/\theta}$, $c>0$, $\theta\in (0,1)$, $i> i_0\ge 0$.
Its modification is Mandelfrot's law (Mandelbrot, 1965) $f_i=c(i+\beta)^{-1/\theta}$, $\beta\ge 0$.  

Probabilistic interpretation of these and similar laws is an infinite urn scheme studied by 
Bahadur (1960), Karlin (1967). There are $n$ balls that are distributed to urns independently and randomly;
there are infinitely many urns. Each ball goes to urn $i$ with probability $p_i>0$, $p_1+p_2+\ldots=1$
(frequences converge a.s. to probabilities). 

We assume that $p_1\geq p_2 \geq \ldots $ and that one of the following asymptotics hold (the second is wider than the first):
\begin{equation}\label{cond2}
p_i=ci^{-1/\theta}(1+o(i^{-1/2})),
\end{equation}
$\theta\in(0,1)$, $c=c(\theta)$ (this assumption includes Zipf's and Mandelbrot's laws);
\begin{equation}\label{cond3}
p_i=i^{-1/\theta}L_0(i,\theta), 
\end{equation}
$\theta\in(0,1)$, $L_0(i,\theta)$ is a slowly varying function of $i$.  
 
Our aim is to construct asymptotically normal estimators of $\theta$ under (\ref{cond2}). We will prove its strong consistency under (\ref{cond3}). 
So we will use statistics that have been studied by Bahadur (1960), Karlin (1967), 
Dutko (1989), Key (1992, 1996), Zakrevskaya $\&$ Kovalevskii (2001), Gnedin, Hansen $\&$ Pitman (2007), 
Boonta $\&$ Neammanee (2007), Hwang $\&$ Janson (2008),
Bogachev, Gnedin $\&$ Yakubovich (2008), Barbour (2009), Barbour $\&$ Gnedin (2009), Ohannessian $\&$ Dahleh (2012), 
Chebunin (2014), Chebunin $\&$ Kovalevskii (2016), 
Muratov $\&$ Zuyev (2016), Ben-Hamou, Boucheron $\&$ Ohannessian (2017).

Let us denote by $J_i(n)$ the number of balls in urn $i$. $R_n$ is the number of nonempty urns, and $R^{*}_{n,k}$ is the number 
of urns with not lesser than $k\geq 1$ balls
\[
 R_n=\sum_{i=1}^\infty{\mathbf I}\{J_i(n)>0\}, \ \ \ R^{*}_{n,k}=\sum_{i=1}^{\infty} {\bf I} (J_i(n)\ge k).
\]
Note that $R^{*}_{n,1}=R_n$. Numbers of urns with exactly $k$ balls: 
$
R_{n,k}=R^{*}_{n,k}-R^{*}_{n,k+1}
$.
The number of urns with odd number of balls 
$$
U_n=\sum_{i=1}^\infty {\bf I}(J_i(n)\equiv 1 ({\rm mod} \ 2)).
$$

Karlin (1967) suggested to study a random sample with a random number of experiments $\Pi(n)$. Here 
$\{ \Pi(t), \ t\ge 0 \}$ is a Poisson process with parameter $1$. Procedure of the random choice of an urn and the Poisson process
are independent.  Processes $\{J_i(\Pi(t))\stackrel{def}{=} \Pi_i(t), \ t\geq 0\}$ are independent Poisson with parameters $p_i$.
Along with the listed papers, the poissonization was used by Ben-Hamou, Boucheron $\&$ Gassiat (2016) in estimating codes on countable alphabets,
by Durieu $\&$ Wang (2016) for proof of functional CLT for some randomization of statistics $R_n$ and $U_n$, 
by Grubel $\&$ Hitczenko (2009) in studying limit distributions of gaps in discrete random samples, 
by Khmaladze (2011) for more general occupancy schemes.

From definition
\[
R^{*}_{\Pi(t),k}=\sum_{i=1}^{\infty} {\bf I}(\Pi_i(t)\geq k), \  R_{\Pi(t),k}=\sum_{i=1}^{\infty} {\bf I}(\Pi_i(t)= k), \ 
U_{\Pi(t)}=\sum_{i=1}^\infty {\bf I}(\Pi_i(t)\equiv 1 ({\rm mod} \ 2)).
\]

Karlin (1967) introduced function
$\alpha(x)=\max\{j|\ p_j\geq 1/x\}$ and proved that (\ref{cond3}) resulted in 
$\alpha(x)=x^{\theta} L(x,\theta)$, 
$L(x,\theta)$ is a slowly varying function as $x \to \infty$. 

Karlin proved SLLNs for all these statistics under (\ref{cond3}). Karlin proved CLTs for $R_n$, $U_n$ and vector 
$(R_{n,1},\ldots, R_{n,d})$ for any finite $d$.

Karlin proved that asymptotics of expectations of all of these statistics is proportional to $\alpha(n)$ with some coefficient depending on $\theta$ only.
This law was found for texts empirically (with $L(x,\theta)=L(\theta)$) by Herdan (1960) and Heaps (1978, Sect. 3.7). 
It is interesting that modern large-scale 
studies of languages demonstrate a deviation from this law (Petersen et al., 2012)  that is interpreted as a decrease of need in new words.

The authors do not know any estimator of $\theta$ with proved asymptotic normality. An estimator of Zakrevskaya $\&$ Kovalevskii (2001) 
founded by a substitution method is (we will see it) asymptotically normal for Zipf's law but authors proved consistency only.
An estimator of Chebunin (2014) is strongly consistent but is not asymptotically normal. We will prove asymptotic normality of 
estimators of Ohannessian $\&$ Dahleh (2012) under (\ref{cond2}) but authors proved only strong consistency under (\ref{cond3}). 

The rest of the paper is organized as follows. In Section 2 we construct asymptotically normal 
estimators of $\theta$ using only one of the statistics. It is possible only if 
we know constant $C$ (it can be a differentiable function of $\theta$) in (\ref{cond2}), and all the estimators in this case are implicit. 
In Section 3 we prove asymptotic normality of estimators that based on two statistics. We use multidimensional CLTs for $(R_{n,1},\ldots, R_{n,d})$ 
that have proved by Karlin (1967) and  for $(R_n,  R_{n,1},\ldots, R_{n,d})$ that we prove in Appendix in a functional generalization.

We use designation $\Rightarrow {\bf N}_{0,\sigma^2}$ for weak convergence to a normal distribution with zero mean and variance $\sigma^2$. 
All convergencies are under $n \to \infty$. 

\section{Implicit estimators that use only one statistics}

We will prove a general theorem for some abstract statistics $S_n$ in the infinite urn scheme 
with neccesary properties Then we will prove these properties to be held for all
statistics under consideration if one assume (\ref{cond2}).
 
 Let $S_n/ n^\theta l(n,\theta)\stackrel{a.s.}{\to} 1$ as $n \to \infty$, where $l(\theta,n)$ is a slowly varying function.
 Let us define $\theta_n^*\in(0,1)$ as a solution of equation 
\begin{equation}\label{p_n}
S_n=n^\theta l(\theta,n).
\end{equation}

As \  
$\ln S_n -\theta \ln n - \ln l(\theta,n) \to 0$, so 
\[
\frac{\ln S_n}{\ln n}   \stackrel{a.s.}{\to} \theta, \ \ \ \textrm{and} \ \ \  
\frac{\ln S_n}{\ln n}-\theta_n^*=\frac{\ln l(\theta_n^*,n)}{\ln n} \stackrel{a.s.}{\to} 0.
\] 
So $\theta_n^*$ is a strongly consistent estimator of $\theta$. 
We  will study asymptotic normality of $\theta_n^*$.
Let 
\begin{equation}\label{S_n}
{\bf E} S_n=n^\theta l(\theta,n)+o(\sqrt{{\bf E} S_n}), \ \ \frac{{\bf Var} S_n}{{\bf E} S_n}\to \sigma^2, \ \   
\frac{S_n}{{\bf E} S_n}\stackrel{a.s.}{\to}1,  \ \  \frac{S_n-{\bf E} S_n}{\sqrt{{\bf Var} S_n}}\Rightarrow {\bf N}_{0,1},
\end{equation}
 $l(\theta,n)$ is a slowly varying function as $n \to \infty$.

\begin{Theorem}\label{Th1}
Let (\ref{S_n}) be held and
$$
 \frac{\ln l(\theta^*_n,n)- \ln l(\theta,n)}{(\theta^*_n-\theta)\ln n}\stackrel{def}{=}\widetilde{l}_n  \stackrel{p}{\to} 0,
$$
$\theta^*_n$ be a solution of (\ref{p_n}).
Then 
$$
\ln n \sqrt{S_n}(\theta^*_n -\theta)\Rightarrow {\bf N}_{0,\sigma^2}.
$$
\end{Theorem}

{\it Proof.} \
$
S_n^0:=\frac{S_n -  n^\theta l(\theta,n)}{\sqrt{{\bf Var} S_n}}\Rightarrow {\bf N}_{0,1}
$.
From (\ref{S_n})
\[
\ln S_n - \ln ( n^\theta l(\theta,n))= \ln\left(1+ \frac{S_n}{ n^\theta l(\theta,n)} -1\right)\stackrel{a.s.}{\sim} \frac{S_n}{ n^\theta l(\theta,n)} -1
\]
 as $n \to \infty$. Then
\[
S_n^0= \frac {n^\theta l(\theta,n)}{\sqrt{{\bf Var} S_n}} \left( \frac{S_n}{n^\theta l(\theta,n)} -1 \right)
\stackrel{a.s.}{\sim} \sqrt{\frac{n^\theta l(\theta,n)}{\sigma^2}}\left( \frac{S_n}{ n^\theta l(\theta,n)} - 1\right)
\]
 \[
\stackrel{a.s.}{\sim} \sqrt{\frac{S_n}{\sigma^2}}\left(\ln S_n - \theta \ln n - \ln l(\theta,n)\right)=
\sqrt{\frac{S_n}{\sigma^2}}\left(\theta^*_n\ln n+\ln l(\theta^*_n,n) - \theta \ln n - \ln l(\theta,n)\right)
\]
\[
=\ln n \sqrt{\frac{S_n}{\sigma^2}}(\theta^*_n -\theta)\left(1+
\frac{\ln l(\theta^*_n, n)  - \ln l(\theta,n)}{(\theta^*_n-\theta)\ln n}  \right)\sim \ln n \sqrt{\frac{S_n}{\sigma^2}}(\theta^*_n -\theta)
\]
in probability as $n\to\infty$.
The theorem is proved.

If $l(\theta, x)=l(\theta)$ is differentiable on $\theta$ then  $\widetilde{l}_n\stackrel{a.s.}{\to} 0$ as $n\to\infty$. Really,
$\theta^*_n \stackrel{a.s.}{\to} \theta$, and  
\[
\widetilde{l}_n= \frac{\ln l(\theta^*_n)- \ln l(\theta)}{(\theta^*_n-\theta)\ln n}
\stackrel{a.s.}{\sim} \frac{l'_\theta(\theta)}{l(\theta) \ln n}\stackrel{a.s.}{\to} 0.
\]

Let $\theta\in(0,1)$, (\ref{cond3}) holds and $L_0(n,\theta) \to c(\theta)$ as $n \to \infty$. Then
$\alpha(x)=\alpha(x,\theta)\sim x^{\theta} c^\theta$.
For example,
\[
p_i(\theta)=\frac{(i-i_0)^{-1/\theta}}{\zeta(1/\theta)}, \ i>i_0,
\]
$i_0$ is integer, $\zeta(z)=\sum_{j=1}^{\infty} j^{-z}$ is Riemann function.
In this case $\alpha(\theta,n)=[(n\zeta(1/\theta))^{\theta}]+i_0$.
From SLLN
$$
\ln R_n - \theta \ln n  -\ln(\Gamma (1-\theta)c^\theta)=\ln n \left(\frac{\ln R_n}{\ln n}- \theta\right) - \ln(\Gamma (1-\theta)c^\theta)\stackrel{a.s.}{\to} 0.
$$
If we use estimator $\theta^*=\ln R_n/\ln n$ (it is consistent, Chebinin (2014)) then $\ln n (\theta^* - \theta)$
goes to a some constant a.s. So we need in implicit estimators for asymptotic normality.
We will base implicit estimators on $R_n$,  $U_n$ or $R_{n,k}$.
Karlin (1967) proved
\[
{\bf E} R_n \sim  \Gamma(1-\theta)  c^\theta  n^{\theta}, \ \ {\bf Var} R_n \sim \left(2^{\theta}-1\right)\Gamma(1-\theta) c^\theta n^{\theta}, \ \ 
\frac{{\bf Var} R_n}{{\bf E} R_n}\to 2^\theta-1,
\]
\[
{\bf E} U_n \sim  2^{\theta-1} \Gamma(1-\theta)  c^\theta n^{\theta}, \ \ 
{\bf Var} U_n\sim 4^{\theta-1} \Gamma(1-\theta)  c^\theta n^{\theta} , \ \   \frac{{\bf Var} U_n}{{\bf E} U_n}\to 2^{\theta-1},
\]
\[
{\bf E} R_{n,k} \sim \theta \frac{\Gamma(k-\theta)}{k!}  c^\theta n^{\theta}, \ \
 {\bf Var} R_{n,k} \sim \frac{\theta}{k!}\left( \Gamma(k-\theta) - \frac{2^\theta\Gamma(2k-\theta)}{2^{2k}k!}\right) c^\theta n^\theta,
\]
\[
\frac{{\bf Var} R_{n,k}}{{\bf E} R_{n,k}}\to 1- \frac{2^\theta\Gamma(2k-\theta)}{2^{2k}k!\Gamma(k-\theta)}.
\]

\begin{Lemma}
If $ \alpha(x)=(c x)^\theta+o(x^\frac\theta2)$ then
\[
{\bf E} R_n =  \Gamma(1-\theta) c^\theta n^{\theta} + o(n^\frac\theta2), \ \
{\bf E} U_n =  2^{\theta-1} \Gamma(1-\theta) c^\theta n^{\theta} +o(n^\frac\theta2), 
\]
\[
{\bf E} R_{n,k} = \theta \frac{\Gamma(k-\theta)}{k!}  c^\theta n^{\theta} +o(n^\frac\theta2).
\]
\end{Lemma}
{\it Proof.} There are convergencies (see Karlin (1967) and Gnedin, Hansen $\&$ Pitman (2007), Lemma 1)
\[
{\bf E} (R_n-R_{\Pi(n)})\to0, \ \ {\bf E} (U_n-U_{\Pi(n)})\to0, \ \ {\bf E} (R_{n,k}-R_{\Pi(n),k})\to0.
\]
We use Karlin (1967) representation, integration by parts and substitution $n t=x$:
\[
{\bf E} R_{\Pi(n)}  
= \int_0^{\infty} \left( 1- e^{-n/x}\right) d\alpha(x)
= \int_0^{\infty}\alpha(x) n x^{-2}e^{- n/x} dx
\]
\[
=  \int_0^{\infty}((c nt)^\theta+o((nt)^\frac\theta2))  t^{-2}e^{- 1/t} dt
=\Gamma(1-\theta) c^\theta n^{\theta} + o(n^\frac\theta2).
\]

Analogously for ${\bf E} U_{\Pi(n)}$ and ${\bf E} R_{\Pi(n),k}$.
Proof is complete.

\begin{Lemma}
If (\ref{cond2}) holds then   
$\alpha(x)=(c x)^\theta+o(x^\frac\theta2)$.
\end{Lemma}
{\it Proof.}
Let us solve equation
$
c\cdot i^{-1/\theta}(1+o(i^{-\frac12}))=\frac1x
$ 
for large enough $x$. 
\[
i=(c x)^\theta (1+o(i^{-\frac12}))^\theta=(c x)^\theta (1+o(i^{-\frac12}))
=(c x)^\theta (1+o((c x)^{-\frac\theta2} (1+o(i^{-\frac12}))^{-\frac12}))
\] 
\[
=(c x)^\theta (1+o((c x)^{-\frac\theta2} (1+o(i^{-\frac12}))))= (c x)^\theta + o(x^\frac\theta2).
\] 
Proof is complete.

\begin{Corollary}
If (\ref{cond2}) holds, $c$ is known, $\frac{dc}{d{\theta}}$ exists, $\theta_{n,R}^*$, $\theta_{n,U}^*$, $\theta_{n,k}^*$ are solutions of equations 
\[ R_n=\Gamma(1-\theta)(cn)^\theta, \ \  U_n= 2^{\theta-1}\Gamma(1-\theta)(cn)^\theta, \ \ 
R_{n,k}=\theta \frac{\Gamma(k-\theta)}{k!} (cn)^\theta
\]
 respectively, then 
\[
\ln n \sqrt{R_n}(\theta^*_{n,R} -\theta)\Rightarrow {\bf N}_{0,2^\theta-1},
\ \
\ln n \sqrt{U_n}(\theta^*_{n,U} -\theta)\Rightarrow {\bf N}_{0,2^{\theta-1}},
\]
\[
\ln n \sqrt{R_{n,k}}(\theta^*_{n,k} -\theta)\Rightarrow {\bf N}_{0,\sigma^2}, \ \  
\sigma^2=1- \frac{2^\theta\Gamma(2k-\theta)}{2^{2k}k!\Gamma(k-\theta)}.
\]
\end{Corollary}

\section{Explicit estimators on a base of two statistics}

Let parameter (function) $c$ be unknown. In this case we need in two statistics to estimate $\theta$.
Some of the following estimators are proposed by Ohannessian $\&$ Dahleh (2012). We will prove its 
asymptotical normality. Note that rates of convergence will be slower in this case.

\begin{Theorem}\label{Th2}
 If $\frac{{\bf E} R_{n,1}-\theta {\bf E} R_n}{\sqrt{\alpha(n)}}\to 0$ then
$
\sqrt{R_n} \left(\frac{R_{n,1}}{R_n} - \theta\right) \Rightarrow {\bf N}_{0,\sigma_0^2},
$
\[
\sigma_0^2=\theta((9\theta-1)2^{\theta-2}+1-\theta).
\]
\end{Theorem}

{\it Proof.}
Using SLLN we have
\[
\sqrt{R_n} \left(\frac{R_{n,1}}{R_n} - \theta\right)= \frac{R_{n,1}-\theta R_n}{\sqrt{R_n}} 
\stackrel{a.s.}{\sim} \frac{R_{n,1}-\theta R_n}{\sqrt{\Gamma(1-\theta) \alpha(n)}}
\]
\[
\stackrel{a.s.}{\sim} \frac{R_{n,1}-{\bf E} R_{n,1} -\theta(R_n-{\bf E} R_{n})}{\sqrt{\Gamma(1-\theta) \alpha(n)}}=
\frac{1}{\sqrt{\Gamma(1-\theta)}} \left(
\frac{R_{n,1}- {\bf E}R_{n,1}}{\sqrt{\alpha(n)}}-\theta
\frac{R_{n}- {\bf E}R_{n}}{\sqrt{\alpha(n)}}\right).
\]
Then we calculate limiting variance using Corollary \ref{cor3}.
Proof is complete.

Note that  $\sigma_0^2<4$  for $\theta\in(0,1)$.

\begin{Theorem}\label{Th3}
If  $\frac{(k-\theta){\bf E} R_{n,k}-(k+1) {\bf E} R_{n,k+1}}{\sqrt{\alpha(n)}}\to 0$ then 
\[
\sqrt{R_{n,k}} \left(\frac{kR_{n,k}-(k+1) R_{n,k+1}}{R_{n,k}} - \theta\right)
\Rightarrow {\bf N}_{0,\sigma_k^2},
\]
\[
\sigma_k^2
=(k-\theta)(2k+1-\theta)-\frac{(2k-\theta+\theta^2)}{k2^{2k+2-\theta} {\rm B}(k-\theta,k)},
\]
${\rm B}$ is a Beta function.
\end{Theorem}
{\it Proof.} 
Using SLLN we have
\[
\sqrt{R_{n,k}} \left(\frac{kR_{n,k}-(k+1) R_{n,k+1}}{R_{n,k}} - \theta\right)=
 \frac{(k-\theta)R_{n,k}-(k+1) R_{n,k+1}}{\sqrt{R_{n,k}}} 
\]
\[
\stackrel{a.s.}{\sim} \frac{(k-\theta)(R_{n,k}-{\bf E} R_{n,k})-(k+1) (R_{n,k+1}-{\bf E} R_{n,k+1})}{\sqrt{\theta \frac{\Gamma(k-\theta)}{k!}  \alpha(n)}}
\]
\[
=
\frac{1}{\sqrt{\theta \frac{\Gamma(k-\theta)}{k!}}} \left(
(k-\theta) \frac{R_{n,k}- {\bf E}R_{n,k}}{\sqrt{\alpha(n)}}-(k+1)
\frac{R_{n,k+1}- {\bf E}R_{n,k+1}}{\sqrt{\alpha(n)}}\right).
\]
Then we calculate limiting variance on the base of Theorem 5 in Karlin (1967).
Proof is complete.

From Lemma 1 and Lemma 2 we obtain the following corollary.

\begin{Corollary}
Assumptions of Theorem \ref{Th2} and Theorem \ref{Th3} are held under (\ref{cond2}).
\end{Corollary}

\section*{Appendix: Functional Central Limit Theorem}

Let for $t\in [0,1], \ k\geq 1$ 
\[
Y_{n,k}^{*}(t) =\frac{R^{*}_{[nt],k}-{\bf E} R^{*}_{[nt],k}}{(\alpha(n))^{1/2}}, 
\ \ \ \ \ \ \ \ \ \  
 Y_{n,k}(t) =\frac{R_{[nt],k}-{\bf E} R_{[nt],k}}{(\alpha(n))^{1/2}}.
\]
\begin{Theorem}\label{Th4} 
Let us assume that (\ref{cond3}) holds, $\nu\geq 1$ is integer. 
Then random process
$
\left((Y^*_{n,1}(t),  
Y_{n,1}(t),\ldots,
Y_{n,\nu}(t)), 
\ 0 \leq t \leq 1 \right)
$
converges weakly in the uniform metrics in  $D(0,1)$ 
to $\nu+1$-dimensional Gaussian process with zero expectation and covariance function 
$(c_{ij}(\tau,t))_{i,j=0}^{\nu}$,
\[
c_{ij}(\tau,t)=
  \frac{\theta \tau^i (t-\tau)^{j-i} t^{\theta-j} \Gamma(j-\theta)}{i!(j-i)!}  
-
\frac{\theta \tau^i t^j (t+\tau)^{\theta-i-j} \Gamma(i+j-\theta)}{i!j!}  \ \ {\rm for} \ \ 1 \leq i < j, \ \tau\leq t,
\]
\[
c_{ij}(\tau,t)=
-
\frac{\theta \tau^i t^j (t+\tau)^{\theta-i-j} \Gamma(i+j-\theta)}{i!j!}  \ \ {\rm for} \ \ i> j\geq 1, \ \tau\leq t,
\]
\[
c_{ii}(\tau,t)=
  \frac{\theta  t^{\theta} \Gamma(i-\theta)}{i!}  
-
\frac{\theta \tau^i t^i (t+\tau)^{\theta-2i} \Gamma(2i-\theta)}{(i!)^2}  \ \ {\rm for} \ \ i> 0, \ \tau\leq t,
\]
\[
c_{00}(\tau,t)=
\left((t+\tau)^{\theta}-t^{\theta}\right) \Gamma(1-\theta)  \ \ {\rm for} \ \  \tau\leq t,
\]
\[
c_{i0}(\tau,t)=
-
\frac{\theta \tau^i (t+\tau)^{\theta-i} \Gamma(i-\theta)}{i!}  \ \ {\rm for} \ \ i> 0, \ \tau\leq t,
\]
\[
c_{0j}(\tau,t)=
\frac{\theta ((t-\tau)^j t^{\theta-j} - t^j (t+\tau)^{\theta-j}) \Gamma(j-\theta)}{j!}  \ \ {\rm for} \ \ j>0, \ \tau\leq t,
\]
$c_{ji}(t,\tau)=c_{ij}(\tau,t)$.
\end{Theorem}

{\it Proof}. 
We base on Theorem 3 in Chebunin $\&$ Kovalevskii (2016) and use formulas
\[
c_{ij}(\tau,t)=c^{*}_{ij}(\tau,t)-c^{*}_{i+1,j}(\tau,t)-c^{*}_{i,j+1}(\tau,t)
+c^{*}_{i+1,j+1}(\tau,t),
\]
\[
c_{0j}(\tau,t)=c^{*}_{1j}(\tau,t)-c^{*}_{1,j+1}(\tau,t), \ \ 
c_{i0}(\tau,t)=c^{*}_{i1}(\tau,t)-c^{*}_{i+1,1}(\tau,t).
\]
Proof is complete.

The limiting $\nu$-dimensional Gaussian process is self-similar with Hurst parameter $H=\theta/2<1/2$. 
Its first component coinsides in distribution with the first component of the limiting process in Theorem~1 in Durieu $\&$ Wang (2015).  

We need in a some specific corollary to calculate limiting variance in Theorem \ref{Th2}.

\begin{Corollary}\label{cor3}
In assumptions of Theorem \ref{Th4}, random vector $\left((Y^*_{n,1}(1), Y_{n,1}(1) \right)$
converges weakly to a normal one with zero mean and covariance matrix
\[
\Gamma(1-\theta)
\left(
\begin{array}{cc}
2^{\theta}-1 & -\theta 2^{\theta-1}\\
-\theta 2^{\theta-1} & \theta(1-2^{\theta-2}(1-\theta))  
\end{array}
\right).
\]
\end{Corollary}

{\bf Acknowledgement}
Our research was supported by RFBR grant 17-01-00683. 


\bigskip

\footnotesize

{\sc Bahadur, R. R.}, 1960. On the number of distinct values in a large sample from an infinite discrete distribution.
Proceedings of the National Institute of Sciences of India, 26A, Supp. II, 67--75.

{\sc Barbour, A. D.}, 2009. Univariate approximations in the infinite occupancy
scheme. Alea 6, 415--433.

{\sc Barbour, A. D.,  Gnedin, A. V.}, 2009.
Small counts in the infinite occupancy scheme.
Electronic Journal of Probability, 
Vol. 14, Paper no. 13, 365--384.

{\sc Ben-Hamou, A., Boucheron, S., Gassiat, E.}, 2016.
Pattern coding meets censoring: (almost) adaptive coding on countable alphabets. 
Preprint. arXiv:1608.08367.

{\sc Ben-Hamou, A., Boucheron, S., Ohannessian, M. I.}, 2017.
Concentration inequalities in the infinite urn scheme for occupancy counts and the missing mass, with applications.
Bernoulli, V. 23, 249--287.

{\sc Bogachev, L. V., Gnedin, A. V., Yakubovich, Y. V.}, 2008.
 On the variance of the number of occupied boxes.
 Adv. Appl. Math., V. 40, 401--432.

{\sc Boonta, S., Neammanee, K.}, 2007. Bounds on random infinite urn model.
Bulletin of the Malaysian Mathematical Sciences Society. Second Series, V. 30.2, 121--128.

{\sc Chebunin, M. G.}, 2014. Estimation of parameters of probabilistic models which is based on the number of different elements 
in a sample. Sib. Zh. Ind. Mat., 17:3, 135--147 (in Russian).

{\sc Chebunin, M., Kovalevskii, A.}, 2016. 
Functional central limit theorems for certain statistics in an infinite urn scheme. 
Statistics and Probability Letters, V. 119, 344--348.

{\sc Durieu, O., Wang, Y.}, 2016. From infinite urn schemes to decompositions of self-similar Gaussian processes.  
Electron. J. Probab., 2016, V. 21, paper No. 43.

{\sc Dutko, M.}, 1989.
Central limit theorems for infinite urn models, Ann. Probab. 17,
1255--1263.

{\sc Gnedin, A., Hansen, B., Pitman, J.}, 2007.
Notes on the occupancy problem with
infinitely many boxes: general
asymptotics and power laws.
Probability Surveys,
Vol. 4, 146--171.

{\sc Grubel, R., and Hitczenko, P.}, 2009.
Gaps in discrete random samples, 
J. Appl. Probab., V. 46, 1038--1051.

{\sc Heaps, H. S.}, 1978.
Information Retrieval: Computational and Theoretical Aspects, Academic Press.

{\sc Herdan, G.}, 1960.
 Type-token mathematics, The Hague: Mouton.

{\sc Hwang, H.-K., Janson, S.}, 2008. Local Limit Theorems for Finite and Infinite 
Urn Models. The Annals of Probability, Vol. 36, No. 3,  992--1022.

{\sc Karlin, S.}, 1967. Central Limit Theorems for Certain Infinite Urn Schemes. 
Jounal of Mathematics and Mechanics, Vol. 17, No. 4,  373--401.

{\sc Key, E. S.}, 1992. Rare Numbers. Journal of Theoretical Probability,
Vol. 5, No. 2, 375--389.

{\sc Key, E. S.}, 1996.
Divergence rates for the number of rare numbers. Journal of Theoretical Probability, Volume
9, No. 2, 413--428.

{\sc Khmaladze, E. V.}, 2011.
Convergence properties in certain occupancy problems including the Karlin-Rouault law,
J. Appl. Probab., V. 48, 1095--1113.

{\sc Mandelbrot, B.}, 1965. 
Information Theory and Psycholinguistics. In B.B. Wolman and E. Nagel. Scientific psychology. Basic Books.

{\sc Muratov, A., and Zuyev, S.}, 2016.
Bit flipping and time to recover,
J. Appl. Probab., V. 53, 650--666.

{\sc Ohannessian, M. I., Dahleh, M. A.}, 2012.
Rare probability estimation under regularly varying heavy tails,   
Proceedings of the 25th Annual Conference on Learning Theory, PMLR 23:21.1--21.24.

{\sc Petersen, A. M., Tenenbaum, J. N.,  Havlin, S., Stanley, H. E., Perc, M.}, 2012. 
Languages cool as they expand: Allometric scaling and the decreasing need for new words.
Scientific Reports 2, Article No. 943.

{\sc Zakrevskaya, N. S.,  Kovalevskii, A. P.}, 2001. One-parameter probabilistic models of text statistics. 
Sib. Zh. Ind. Mat., 4:2, 142--153 (in Russian).

{\sc Zipf, G. K.}, 1949.
Human behavior and the principle of least effort. Cambridge: Univ. Press.

\end{document}